\documentclass[a4paper,11pt]{article}
\usepackage{bbm,amssymb,amsthm,amsmath,amsfonts,fullpage,algorithm,algorithmic}
\usepackage{color,epsfig}
\usepackage{enumerate}
\def\R{\mathbb{R}}

\def\N{\mathbb{N}}

\def\E{\mathbb{E}}

\def\1{\mathbbm{1}}

\DeclareMathOperator*{\argmax}{\arg\,\max}
\DeclareMathOperator*{\supp}{supp}

\def\X{\mathcal{X}}

\def\uq{\underline{q}}
\def\G{\mathcal{G}}

\def\taut{\tilde{\tau}}
\def\Tt{\tilde{T}}
\def\Ft{\tilde{F}}
\def\Rh{\hat{R}}

\newtheorem{lemma}{Lemma}
\newtheorem{corollary}{Corollary}
\newtheorem{proposition}{Proposition}
\newtheorem{theorem}{Theorem}
\newtheorem{remark}{Remark}

\title{Optimal discovery with probabilistic expert advice}

\author{S{\'e}bastien Bubeck \\
Department of Operations and Financial Engineering\\
Princeton University \\
Princeton, NJ, 08544, USA \\
{\tt sbubeck@princeton.edu }
\\ \\
 Damien Ernst \\
 Universit\'e of Liege \\ 
Liegem B-4000, Belgium \\
{\tt dernst@ulg.ac.be}
\\ \\
Aur\'elien Garivier \\
LTCI, CNRS \& Telecom ParisTech \\ 
Paris, France \\
{\tt garivier@telecom-paristech.fr }}

\begin{document}

\maketitle

\begin{abstract}
We consider  an original problem that arises from the issue of security analysis of a power system and that we name  optimal discovery with probabilistic expert advice. We address it with an algorithm based on the optimistic paradigm and the Good-Turing missing mass estimator. We show that this strategy uniformly attains the optimal discovery rate in a macroscopic limit sense, under some assumptions on the probabilistic experts. We also provide numerical experiments suggesting that this optimal behavior may still hold under weaker assumptions.

\end{abstract}

\section{Introduction}
In this paper we consider the following problem: Let $\X$ be a set, and $A \subset \X$ be a set of interesting elements in $\X$. One can access $\X$ only through requests to a finite set of probabilistic experts. More precisely, when one makes a request to the $i^{th}$ expert, the latter draws independently at random a point from a fixed probability distribution $P_i$ over $\X$. One is interested in discovering rapidly as many elements of $A$ as possible, by making sequential requests to the experts.

\subsection{Motivation}
The original motivation for this problem arises from the issue of real-time security analysis of a power system. This problem often amounts to identifying in a  set of  ‘credible’ contingencies those that may indeed endanger the security of the power system and perhaps lead to a  system collapse with catastrophic consequences (e.g., an entire region, country may be without electrical power for hours). Once those dangerous contingencies have been identified,   the system operators usually take preventive actions so as to ensure that they could mitigate their effect on the system in the likelihood they would occur. Note that usually, the dangerous contingencies are very rare with respect to the non dangerous ones.
A straightforward approach for tackling this security analysis problem is to simulate the power system dynamics for every credible contingency so as to identify those that are indeed dangerous.
Unfortunately, when the set of credible contingencies contains a large number of elements (say, there are more than $10^5$ ‘credible’ contingencies) such an approach may not possible anymore since the computational resources required to simulate every contingency may excess those that are usually available during the few (tens of) minutes available for the real-time security analysis. One is therefore left with the problem of identifying within this short time-frame a maximum number of dangerous contingencies rather than all of them.
The approach proposed in \cite{Fonteneau11phd,FonteneauErnstDruetPanciaticiWehenkel10power} addresses this problem by building first very rapidly what could be described as a probability distribution $P$ over the set of credible contingencies that points with  significant probability to contingencies which are dangerous. Afterwards, this probability distribution is used to draw the contingencies to be analyzed through simulations. When the computational resources are exhausted, the approach outputs the contingencies found to be dangerous.
One of the main shortcoming of this approach is that usually $P$ points only with a significant probability to a few of the dangerous contingencies and not all of them. This in turn makes this probability distribution not more likely to generate after a few draws new dangerous contingencies than for example a uniform one. The dangerous contingencies to which $P$ points to with a significant probability depend however strongly on the set of (sometimes arbitrary) engineering choices that have been made for building it. One possible strategy to ensure that more dangerous contingencies can be identified within a limited budget of draws would therefore be to consider $K>1$ sets of engineering choices  to build $K$ different probability distributions $P_1$, $P_2$, $\ldots$, $P_K$ and to draw the contingencies from these $K$ distributions rather than only from a single one. This strategy raises however an important question to which this paper tries to answer: how should the distributions be selected  for being able to generate with a given number of draws a maximum number of dangerous contingencies?
We consider the specific case where the contingencies are sequentially drawn and where the distribution selected for generating a contingency at one instant can be based on the past distributions that have been selected, the contingencies that have been already drawn and the results of the security analyses (dangerous/non dangerous) for these contingencies.
%This problem is carefully formalized in the next subsection of this paper.
This corresponds exactly to the optimal discovery problem with expert advice described above.
We believe that this framework has many other possible applications, such as for example web-based content access.

\subsection{Setting and notation} \label{sec:settingnotation}
In this paper we restrict our attention to finite or countably infinite  sets $\X$. %, and probabilistic experts with non-intersecting supports. 
We denote by $K$ the number of experts. For each $i\in\{1,\dots,K\}$, we assume that $(X_{i,n})_{n\geq 1}$ are random variables with distribution $P_i$ such that the $(X_{i,n})_{i,n}$ are independent. Sequential discovery with probabilistic expert advice can be described as follows: at each time step $t\in\N^*$, one picks an index $I_t\in\{1,\dots, K\}$, and one observes $X_{I_t, n_{I_t,t}}$, where 
\[n_{i,t} = \sum_{s\leq t} \1_{\{I_s = i\}}\;.\]
The goal is to choose the $(I_t)_{t\geq 1}$ so as to observe as many elements of $A$ as possible in a fixed horizon $t$, or equivalently to observe all the elements of $A$ within as few time steps as possible. The index $I_{t+1}$ may be chosen according to past observations: it is a (possibly randomized) function of $(I_1, X_{1,I_1}, \dots, I_{t}, X_{I_t, n_{I_t,t}})$. We are mainly interested in the number of interesting items found by the strategy after $t$ time steps:
$$\sum_{x \in A} \1\bigg\{x \in \{X_{1,1},\hdots,X_{1,n_{1,t}}, \hdots, X_{K,1},\hdots,X_{K,n_{K,t}}\}\bigg\}.$$
Note in particular that it of no interest to observe twice the same same element of $A$.
%, and the time $\bar{\tau}(f)$ 
%at which the strategy found the $f^{th}$ interesting item:
%$$\bar{\tau}(f) = \min\{t\geq 0 : \bar{F}(t) = f\} .$$

While Algorithm Good-UCB, presented in Section~\ref{sec:goodUCB}, can be used in a more general setting (as illustrated in Section~\ref{sec:simus}), for the mathematical analysis we restrict our attention to the case of probabilistic experts with the following properties: 
\begin{enumerate}[(i)]
\item non-intersecting supports: $A \cap \supp(P_i) \cap \supp(P_j) = \emptyset$ for $i \neq j$,
\item finite supports with the same cardinality: $|\supp(P_i)| = N, \forall i \in \{1,\hdots,K\}$,
\item uniform distributions: $P_i(x) = \frac{1}{N}, \forall x \in \supp(P_i), \forall i \in \{1,\hdots,K\}$.
\end{enumerate}
These asumptions are made in order to be able to compare the performance of the Good-UCB algorithm to an ``oracle'', described below.
Indeed, in that case, this oracle has a very simple behavior.
In this setting it is convenient to reparametrize slightly the problem (in particular we make explicit the dependency on $N$ for reasons that will appear later). Let $\X^N = \{1,\dots,K\} \times \{1,\hdots,N\}$, $A^N\subset \X^N$ the set of interesting items of $\X^N$, and $Q^N = |A^N|$ the number of interesting items. We assume that, for expert $i \in \{1,\dots,K\}$,  $P^N_i$ is the uniform distribution on $\{i\}\times \{1,\hdots,N\}$. We also denote by $Q_i^N = \left|A^N \cap \left(\{i\}\times \{1,\hdots,N\}\right) \right|$ the number of interesting items accessible through requests to expert $i$.
Further notation is given in Section \ref{sec:oracle}.

\subsection{Contribution and content of the paper}
This paper contains the description of a generic algorithm for the optimal discovery problem with probabilistic expert advice, and a theoretical proof of optimality in a particular setting. In Section~\ref{sec:goodUCB}, we first depict our strategy, termed Good-UCB. This algorithm relies on the \emph{optimistic paradigm}  (which led to the UCB (Upper Confidence Bound) algorithm for multi-armed bandits, \cite{AuerEtAl02FiniteTime}), and on a finite-time analysis of the Good-Turing estimator for the missing mass. In order to analyze and quantify the performance of this strategy, we compare it with the \emph{oracle} (closed-loop) policy, a virtual algorithm that would be aware, at each time, of the probability of each item under each distribution, and would thus be able to sample optimally. This strategy is carefully described  and analyzed in Section~\ref{sec:oracle}. The analysis is performed under the non-intersecting and uniform draws assumptions [(i), (ii), (iii)] described above, and in a macroscopic limit sense, that is when the size of the set $\X$ grows to infinity while maintaining a constant proportion of interesting items. More precisely we prove the following theorem, where $F^N(t)$ is the number of interesting items found by the oracle policy after $t$ time steps. %(respectively the time at which the oracle policy found the $f^{th}$ interesting item).

\begin{theorem}\label{th:convUnifOOL}
Assume that, for all $i\in\{1,\dots, K\}, Q^N_i/N$ converges to $q_i\in]0,1[$ as $N$ goes to infinity. 
Then, almost surely, the sequence of mappings $t\mapsto F^N\left([Nt]\right)/N$ converges uniformly on $\R_+$ to a limit denoted $F$ as $N$ goes to infinity.
\end{theorem}

In Section~\ref{sec:oracle} we also give an explicit expression for the limit $F$. Section~\ref{sec:OOL} presents a study of the oracle \emph{open-loop} policy which is defined as the optimal fixed allocation. In this problem it turns out that the oracle open-loop policy achieves the same performance as the oracle closed-loop policy, which in turns yields another formula for the macroscopic discovery rate $F$. In particular these formulas allow to see easily the difference in the macroscopic behavior between optimal policies, and suboptimal policies such as uniform requests, see Remark \ref{rmk:1} for more details. 
\newline

The main result of the paper is given in Section~\ref{sec:goodOptimal}. We show that Good-UCB is a macroscopically optimal policy, that is, the performances of Good-UCB tends to the performances of the oracle policy. More precisely let $\Ft^N(t)$ be the number of interesting items found by Good-UCB after $t$ time steps. % (respectively the time at which Good-UCB found the $f^{th}$ interesting item).
\begin{theorem}\label{th:convUnifGUCB}
Assume that, for all $i\in\{1,\dots, K\}, Q^N_i/N$ converges to $q_i\in]0,1[$ as $N$ goes to infinity. 
Then, almost surely, the sequence of mappings $t\mapsto \Ft^N\left([Nt]\right)/N$ converges uniformly on $\R_+$ to the limiting proportion $F$ found during the same time by the oracle policy.
\end{theorem}

Section~\ref{sec:simus} reports  experimental results that  show that
the Good-UCB  algorithm performs very  well, even in a  setting where
assumptions (i), (ii) and  (iii) are not satisfied anymore. Finally,
Section \ref{section:conclusions} concludes.

\section{The Good-UCB algorithm}~\label{sec:goodUCB}
We describe here the Good-UCB strategy. This algorithm is a sequential method estimating at time $t$, for each expert $i\in\{1,\dots,K\}$, the total probability of the interesting items that remain to be discovered through requests to expert $i$. This estimation is done by adapting the so-called Good-Turing estimator for the missing mass. Then, instead of simply using the distribution with highest estimated missing mass, which proves hazardous, we make use of the \emph{optimistic paradigm} (see  \cite{Agrawal:95,AuerEtAl02FiniteTime} and references therein), a heuristic principle well-known in reinforcement learning, which entails to prefer using an \emph{upper-confidence bound} (UCB) of the missing mass instead. At a given time step, the Good-UCB algorithm simply makes a request to the expert with highest upper-confidence bound on the missing mass at this time step.
We start with the Good-Turing estimator and a brief study of its concentration properties. Then we describe precisely the Good-UCB strategy.

\subsection{Estimating the missing mass}
Our algorithm relies on an estimation at each step of the probability of obtaining a new interesting item by making a request to a given expert. A similar issue was addressed by I.~Good and A.~Turing  as part of their efforts to crack German ciphers for the Enigma machine during World War II. In this subsection, we describe a version of the Good-Turing estimator adapted to our problem.
Let $\Omega$ be a discrete set, and let $A$ be a subset of interesting elements of $\Omega$.
Assume that $X_1,\dots,X_n$ are elements of $\Omega$ drawn independently under the same distribution $P$, and define for every $x\in\Omega$:
\[O_n(x) = \sum_{m=1}^n \1\{X_m=x\},\quad Z_n(x) = \1\{O_n(x) = 0\},\quad U_n(x) = \1\{O_n(x) = 1\}\;. \]
Let $p_{\max} = \max\{P(x) : x\in\Omega\}$, let $R_n = \sum_{x\in A} Z_n(x)P(x)$ denote the missing mass of the interesting items, and let $U_n = \sum_{x\in A} U_n(x)$ the number of elements of $A$ that have been seen exactly once (in linguistics, they are often called \emph{appaxes}).
The idea of the Good-Turing estimator (\cite{Good53GT}, see also \cite{McallesterSchapire00GoodTuring,orlitsky2003always} and references therein) is to estimate the (random) ``missing mass'' $R_n$, which is the total probability of all the interesting items that do not occur in the sample $X_1,\dots,X_n$, by the ``fraction of appaxes'' $\Rh_n = U_n/n$.
This estimator is well-known in linguistics, for instance in order to estimate the number of words in some language, see~\cite{Gale95GT}.
For our particular needs, we derive (using similar techniques as in \cite{McallesterSchapire00GoodTuring}) the following upper-bound on the estimation error:
\begin{proposition}\label{prop:inegGT}
With probability at least $1-\delta$,
 \[ \Rh_n -\frac{1}{n} - \sqrt{\frac{(2/n+p_{\max})^2n\log(2/\delta)}{2}} \leq R_n  \leq  \Rh_n+ \sqrt{\frac{(2/n+p_{\max})^2n\log(2/\delta)}{2}} \]
\end{proposition}

\textbf{Proof:} The random variable  $W_n = R_n-\Rh_n$ is a function of the independent observations $X_1,\dots,X_n$ such that, denoting $W_n = f(X_1,\dots,X_n)$, modifying just one observation has limited impact: $\forall l\in\{1,\dots,n\}, \forall (x_1,\dots,x_n,x'_l)\in\Omega^{n+1}$,
\[
\left|f(x_1,\dots,x_n) - f(x_1,\dots,x_{l-1},x'_l,x_{l+1}, \dots,x_n)\right| \leq \frac{2}{n}+ p_{\max}
\]
By applying McDiarmid's inequality~\cite{Mcdiarmid89boundedDiff}, one gets that, with probability at least $1-\delta$, 
\[\left|W_n - \E[W_n]\right| \leq \sqrt{\frac{(2/n+p_{\max})^2n\log(2/\delta)}{2}}\;.\]
Moreover, 
\begin{align*}
 \E[W_n] &= \sum_{x\in A} \left[P(x)\left(1-P(x)\right)^n - \frac{1}{n}\times nP(x)\left(1-P(x)\right)^{n-1}\right]\\
 & = -\frac{1}{n} \sum_{x\in A} P(x) \times nP(x)\left(1-P(x)\right)^{n-1} \\
 %&= -\frac{1}{n} \sum_{x\in A}P(x) \E\left[ U_n(x)  \right] \\
 & = -\frac{1}{n}\E\left[ \sum_{x\in A}P(x) U_n(x)\right]\in \left[-\frac{1}{n},0\right]\;,
\end{align*}
which concludes the proof.

%Note that when $P(A)$ is small, a significant improvement can be obtained by using Boucheron, Massart and Lugosi's inequality~\cite{BoucheronLugosiMassart09selfBounding}: in fact, if $X'_l$ is drawn from $P$  independently of the $X_l$,
%\[
%\E\left[\left|f(X_1,\dots,X_n) - f(X_1,\dots,X_{l-1},X'_l,X_{l+1},X_n)\right|^2\right] \leq \left(\frac{2}{n}+ p_{\max}\right)^2 P(A)\,
%\]
%and hence one gets that 
%\[\left|W_n - \E[W_n]\right| \leq \sqrt{4(2/n+p_{\max})^2P(A)n\log(2/\delta)}\;.\]

\subsection{The Good-UCB algorithm}\label{sec:GoodTuringEstimator}
Following the example of the well-known Upper-Confidence Bound procedure for multi-armed bandit problems, we propose Algorithm~\ref{algo:GUCB}, which we call Good-UCB in reference to the estimated procedure it relies on.
For every arm $i\in\{1,\dots,K\}$ and for every $t\in\N$, denote 

$O_{i,t}(x) = \sum_{s=1}^{n_{i,t}} \1\{X_{i,s}=x\}$, $\quad O_{t}(x) = \sum_{i=1}^K O_{i,t}(x)$, $\quad U_{i,t}(x) = \1\{O_{i,t}(x) = O_{t}(x) = 1\}$, $\quad U_{i,t} = \sum_{x \in A} U_{i,t}(x).$

For each arm $i\in\{1,\dots,K\}$, the index at time $t$ is composed of the estimate
\[\Rh_{i, t-1}=\frac{U_{i, t-1}}{n_{i, t-1}}\]
 of the missing mass 
\[\sum_{x\in A \setminus \{X_{I_{1}, n_{I_1, 1}}, \dots, X_{I_{t-1}, n_{I_{t-1}, t-1}}\}} P^N_i(x)
\] 
 inflated by a confidence bonus of order $\sqrt{\log(t)/n_{i,t-1}}$.
Good-UCB relies on a tuning parameter $c$ which is discussed below.

\begin{algorithm}
\caption{Good-UCB}
\label{algo:GUCB}
\begin{algorithmic}[1]
\STATE For $1\leq t\leq K$ choose $I_t = t$.
\FOR {$t \geq K+1$}
\STATE Choose $I_t =\argmax_{1\leq i\leq K} \left\{\Rh_{i_,t-1} + c\sqrt{\frac{\log{(t)}}{n_{i,t-1}}} \right\}$
\STATE Observe $X_t$ distributed as $P_{I_t}$ and update $\Rh_{1,t}, \hdots, \Rh_{K,t}$ accordingly
\ENDFOR
\end{algorithmic}
\end{algorithm}

Note that the Good-UCB algorithm is designed for more general probabilistic experts than those satisfying assumptions [(i), (ii), (iii)]. In particular since we do not make the non-intersecting supports assumption (i), the missing mass of a given expert $i$ depends explicitly on the outcomes of {\em all} requests (and not only requests to expert $i$). Note also that the bounds of Proposition~\ref{prop:inegGT} hold for all discrete distributions. The experiments of Section~\ref{sec:simus} validate these observations, and show that Good-UCB behaves very well even when assumptions [(i), (ii), (iii)] are not met.
%Note that the Good-UCB algorithm does not require any assumption on the probability distributions $(P_i)_i$: in fact, the bounds of %Proposition~\ref{prop:inegGT} hold for all discrete distributions, and experiments show that the Good-UCB algorithm behaves very well %even when assumptions (i), (ii) and (iii) are not met.
However, for the theoretical analysis of our algorithm, we focus on large values of $N$ under the non-intersecting and uniform draws  assumptions [(i), (ii), (iii)]: indeed in that case the performance of the oracle strategy is simple and deterministic, so that the optimality of the Good-UCB algorithm can be analyzed.
More precisely, we will show that, in the macroscopic limit, the number of items found at each time by Good-UCB converges to the number of items found by the closed-loop oracle strategy that knows the number of interesting items to find with each expert, at every time, and that may use this information to make its choice. In order to prove this, we first analyze the performance of such an oracle strategy.

\section{The closed-loop oracle strategy}\label{sec:oracle}
From now on we restrict our attention to the setting described in Section \ref{sec:settingnotation} with sets $\X^N$, $A^N$, and experts distributions $P_1^N, \hdots, P_K^N$. Denote by $B^N_i$ the set of interesting items supported by $P^N_i$: $B^N_i = \{x\in\{1,\hdots,N\} : (i,x)\in A^N\}$. Let $Q^N_i = |B^N_i|$; in particular, note that $Q^N = Q^N_1+\dots+Q^N_K$. Without loss of generality, we will assume in the analysis that $Q^N_1\geq Q^N_2\geq\dots\geq Q^N_K$. Successive draws of expert $i$ are denoted $(i,X^N_{i,1}),(i,X^N_{i,2}),\dots$, where the variables $(X^N_{i,n})_{i,n}$ are assumed to be independent.
We denote by $(D^N_{i,k})_{1\leq k\leq Q^N_i}$ the increasing sequence of the indices corresponding to draws for which new interesting items are discovered with expert $i$: 
\[D^N_{i,1} = \min\left\{n\geq 1 : X^N_{i,n} \in B^N_i\right\},
\quad D^N_{i,2} = \min\left\{n\geq D^N_{i,1} : X^N_{i,n} \in B^N_i\setminus \left\{X^N_{i,D^N_{i,1}}\right\} \right\}, \dots\]
Reciprocally, we denote $F^N_{i,n} = \max\{k\in\N : D^N_{i,k}\leq n\}$ the number of items found in the first $n$ draws. We also define $S^N_{i,0} = 0$ and for $k\geq 1, S^N_{i,k} = D^N_{i,k}-D^N_{i,k-1}$. The random variables $S^N_{i,k}$ ($1\leq i\leq K, k\geq 1$) are independent with geometric distribution $\G((1+Q^N_i-k)/N)$.
In particular, for all $k\geq 1$,
\begin{equation}\label{eq:espT}
\E\left[D^N_{i,k}\right]  = \frac{N}{Q^N_i} + \dots+\frac{N}{1+Q^N_i-k}\;.
\end{equation}

\subsection{Description of the closed-loop oracle policy}
When the values of $Q^N_1,\dots,Q^N_K$ are known, so that the number of interesting items to find with each expert is known at every step, an horizon-free optimal closed-loop strategy (denoted in the following as the ``oracle closed-loop strategy'' or as OCL) consists in making a request, at each time step, to one of the experts with highest number of still undiscovered interesting items. 
Hence, an OCL strategy can:
\begin{itemize}
 \item first request expert $1$ for $D^N_{1,Q^N_1-Q^N_2}$ steps;
 \item then, alternatively request
 \begin{itemize}
   \item expert $1$ for $S^N_{1,1+Q^N_1-Q^N_2}$ steps;
   \item expert $2$ for $S^N_{2,{}1}$ steps;
   \item expert $1$ for $S^N_{1,2+Q^N_1-Q^N_2}$ steps;
   \item expert $2$ for $S^N_{2,2}$ steps;
   \item and so on, until there are only $Q^N_3$ undiscovered interesting items on experts $1$ and $2$.
 \end{itemize}
 \item and so on, including successively experts $3,4,\dots,K$ in the alternance.
\end{itemize}

For every $l\in\{0,\dots,Q^N_1\}$, we shall be particularly interested in the waiting time $T^N(l)$ until, under the OCL strategy described above, all experts have at most $l$ undiscovered interesting items.
Obviously, 
\[T^N(l) = \sum_{i: Q^N_i>l} D^N_{i, Q^N_i-l}\;,\]
so that, by Equation \eqref{eq:espT},
\begin{equation}\label{eq:espTl}
\E\left[T^N(l)\right]  = \sum_{i: Q^N_i>l} \frac{N}{Q^N_i} + \dots+\frac{N}{l+1} \;. 
\end{equation}
At that time, the number of items discovered so far is $G^N(l) = \sum_{i=1}^K \left( Q^N_i-l \right)_+$.
For every $f \in \{0,\dots, Q^N\}$, we define $L^N(f)$ to be the maximal number of undiscovered items remaining on an expert, once that $f$ interesting items have been discovered altogether. In particular, note that $L^N(G^N(l)) = l$ for all $l$, and that $G^N(L^N(f))\in\{f-K+1,\dots,f\}$ for all $f$.
Besides, the time $\tau^N(f)$ required by the OCL strategy to collect $f$ interesting items satisfies 
\begin{equation}\label{eq:encadrementTau} T^N\left( L^N(f) \right) \leq \tau^N(f) < T^N\left( L^N(f)-1 \right) \;. 
\end{equation}

The performance of the  OCL strategy is maybe more explicitely expressed by the pseudo-inverse of the mapping of $\tau^N$~: for every integer $t$ let $F^N(t)$ denote the total number of items found up to time $t$. $\tau^N$ and $F^N$ are related as follows:
\begin{align}\label{eq:lientauF}
 \forall t\geq 0,\; &F^N(t) = \max\{f\in\{1,\dots,Q^N\} : \tau^N(f)\leq t\}, \nonumber\\
\forall f\in\{0,\hdots, Q^N\},\; &\tau^N(f) = \min\{t\geq 0 : F^N(t) = f\}\;.
\end{align}

\subsection{Macroscopic limit}\label{sec:macrolimit}
We consider a macroscopic limit where $N$ goes to infinity together with the $Q^N_i$ in such a way that $Q^N_i/N \to q_i\in]0,1[$. Let $q = q_1+\dots+q_K$.
We will show that if $f^N$ is a sequence of integers such that  $f^N/N\to\phi\in]0,q[$, then the normalized waiting time $\tau^N(f^N)/N$ converges as $N$ goes to infinity to a deterministic limit $\tau(\phi)$ that we will compute as a function of $\phi$ and $q_1,\dots,q_K$.
We start with some notation.
Let $G:[0,q_1]\to[0,q]$ be the strictly decreasing mapping defined by
\[G(\lambda) = \sum_{i=1}^K (q_i-\lambda)_+\;,\]
and let $L:[0,q]\to[0,q_1]$ be the inverse mapping.
The following lemmas are proved in the appendix. The second lemma is the key step of the macroscopic analysis.
\begin{lemma}\label{lem:convFL}
The mappings $G$ and $L$ are the  limits of the sequences of mappings $(G^N)_N$ and $(L^N)_N$, respectively, in the following sense:
\begin{itemize}
 \item [(i)] If $l^N\in[0,Q^N_1]$ defines a sequence of integers, then  $l^N/N \to  \lambda\in[0,q_1]$ if and only if $G^N(l^N)/N \to G(\lambda)$ as $N$ goes to infinity.
 \item [(ii)] If $f^N\in[0,Q^N]$ defines a sequence of integers, then  $f^N/N\to\phi\in[0,q]$ if and only if $L^N(f^N)/N \to L(\phi)$ as $N$ goes to infinity.
 \end{itemize}
\end{lemma}
\begin{lemma}\label{lem:convT}
For every $\lambda\in]0, q_1]$, let 
\[T(\lambda) = \sum_{i: q_i>\lambda} \log\frac{q_i}{\lambda} \;.\]
For every sequence $(l^N)_N$ such that $l^N/N$ converges to $\lambda$ as $N$ goes to infinity, $T^N(l^N)/N$ converges almost surely to $T(\lambda)$ as $N$ goes to infinity.
\end{lemma}

\begin{theorem}\label{th:convTau}
For every sequence $(f^N)_N$ such that $f^N\in\{0, \dots,Q^N\}$ and $f^N/N$ converges to $\phi$ as $N$ goes to infinity, $\tau^N(f^N)/N$ converges a.s. to $\tau(\phi) = T(L(\phi))$ as $N$ goes to infinity.
\end{theorem}
\textbf{Proof:}
By Equation~\eqref{eq:encadrementTau}, 
\[ \frac{T^N\left( L^N(f^N) \right)} {N} \leq \frac{\tau^N(f^N)} {N} < \frac{T^N\left( L^N(f^N)-1 \right)}{N}\;.\]
By Lemma~\ref{lem:convFL}, $\lim_{N\to\infty} L^N(f^N)/N  = \lim_{N\to\infty} ( L^N(f^N) - 1 )/N  = L(\phi) $ and thus, by Lemma~\ref{lem:convT}, 
\[\lim_{N\to\infty} \frac{T^N\left( L^N(f^N) \right)} {N}   = \lim_{N\to\infty} \frac{T^N\left( L^N(f^N)-1 \right)}{N} = T(L(\phi))\;.\]
\begin{corollary}\label{cor:convF}
For every sequence sequence $t^N$ of integers such that $t^N/N\to t$ as $N$ goes to infinity, $F^N(t^N)$ converges to $F(t)$, where the function $F:\R\to[0, q[$ is $\tau^{-1} = G\circ T^{-1} $, i.e. 
\[ F(t)  = \sum_{i=1}^K \left( q_i-T^{-1}(t) \right)_+\;.\]
\end{corollary}
The proof, very similar to the previous proofs, is omitted. Another expression for $F$ is obtained in Section~\ref{sec:OOL} in the analysis of the \emph{open-loop} oracle policy.
This allows us to finish the proof of Theorem~\ref{th:convUnifOOL}:
%The sequence of processes $t\mapsto F^N([Nt])/N$, defined on $\R^+$, almost surely converges  to the deterministic function $F$ uniformly on every compact.
%\end{cor}

\textbf{Proof of Theorem~\ref{th:convUnifOOL}~:}
The fact that, almost surely, the sequence of increasing processes $t\mapsto F^N([Nt])/N$ uniformly converges to $F$ on every compact of $\R_+$ is a consequence of Corollary~\ref{cor:convF} and Dini's Theorem. This is sufficient, since $F$ is upper-bounded by $q$.
% 
% In fact, the functionnal law of large numbers (for triangular arrays) implies that for all $i\in\{1,\dots,K\}$,
% \[ \sup_{\lambda \in]0,q_1[} \left\{N^{-1} \left[ T^N_{i,Q^N_i-[N\lambda]} - \E\left[ T^N_{i,Q^N_i-[N\lambda]} \right] \right] \right\}\to 0 \;.\]
% Besides, by Equations~\eqref{eq:espTl} and~\eqref{eq:harmonic}
% \[N^{-1}\E\left[ T^N_{i,Q^N_i-[N\lambda]}\right] = \frac{1}{Q^N_1} + \dots+\frac{1}{[N\lambda]+1} = \log\frac{Q^N_i}{[N\lambda]} + O\left( \frac{1}{N} \right) = \log\frac{q_i}{\lambda}  + O\left( \frac{1}{N} \right) \;. \]
% 
% ATTENTION PETIT PROBLEME QUAND LAMBDA TEND VERS 0, CAR LE O(1/N) DEPEND DE LAMBDA : IL FAUT PEUT-ETRE PRENDRE LE SUP SUR $]\epsilon,q_1[$.
% ]
% Hence, by Equation~\eqref{eq:encadrementTau}, it holds that for all $\phi\in]0,q[$,
% \[\sup_{\lambda \in]0,q_1[} \left\{N^{-1}\tau^N\left([N\phi]\right) - \tau(\phi)\right\}\to 0\,\]
% where $\tau(\phi) = T(L(\phi))$ and $L:]0,q[\mapsto]0,q_1[$ denotes the inverse of the one-to-one, decreasing map $F: \lambda \mapsto \sum_{i=1}^K (q_i-\lambda)_+$. 
% It appears that $\tau$ is the strictly increasing mapping from $]0,q[$ to $]0,+\infty[$ that is defined by:XXX

% 
% As a consequence, the inverse process mapping every (rescaled) time to the proportion of items found converges to the inverse of $\tau$, which is exactly XXX :  hence, in the macroscoping limit there is no difference between closed-loop and open-loop policies.

\section{Macroscopic optimality of the Good-UCB algorithm}\label{sec:goodOptimal}
After $n$ requests to expert $i\in\{1,\dots,K\}$, denote by $\underline{R}^N_{i,n} = (Q^N_i - F^N_{i}(n))/N$  the proportion of interesting items not yet found with that expert.
To estimate this number, we use the Good-Turing estimator 
$$\underline{\Rh}^N_{i,n} = \frac{\sum_{x \in B^N_i} \1\{\sum_{m=1}^n \1\{X_{i,m}=x\}=1\}}{n}\; .$$ 
In particular note that under assumption (i), the estimator $\Rh_{i,t}$ defined in Section~\ref{sec:GoodTuringEstimator} satisfies $\Rh_{i,t} = \underline{\Rh}_{i,n_{i,t}}$. To simplify the proof, we use here a slightly different confidence bonus than the one proposed on line 3 of the algorith. We define here the following upper confidence bound:
\[u^N_{i,n} = \underline{\Rh}^N_{i,n} + \sqrt{\frac{(2/n+1/N)^2 \,n \log(2KN^4)}{2}} \; .\]
%on line 3 of the algorithm.
%The exact form of the confidence bonus used in this proof slightly differs from the suggestion made in~Algorithm~\ref{algo:GUCB} for %the simplicity of the proof, yet the order of magnitude is the same.
According to Proposition~\ref{prop:inegGT}, the event $C^N$ such that: 
\begin{equation}\label{eq:defC}
\forall i\in\{1,\dots,K\}, \forall n\in\{1,N^2\},  u^N_{i,n} - 2\sqrt{\frac{(2/n+1/N)^2 \, n \log(2KN^4)}{2}} -\frac{1}{n} \leq  \underline{R}^N_{i,n} \leq u^N_{i,n} \; ,
\end{equation}
has probability at least $1-N^{-2}$.

The Good-Turing optimistic strategy consists in making a request, at each step, to the expert maximizing the upper-confidence bound $u^N_{i,n_{i,t}}$, where $n_{i,t}$ denotes the (random) number of requests to expert $i$ before time $t$.
Denote by $\taut^N(f)$ the time required by this algorithm to collect $f$ interesting items, and by $\Ft^N(t)$ the number of interesting items found by Good-UCB in the first $t$ rounds.

\begin{theorem}\label{th:optGUCB}
 In the macroscopic limit, $\taut^N(f^N)/N$ converges a.s. to the same limit $\tau(\phi)$ as for the oracle closed-loop policy when $f^N/N$ tends to $\phi$.
Moreover, $\Ft^N(t^N)/N$ converges a.s. to the limiting proportion $F(t)$ found during the same time by the OCL policy when $t^N/N$ converges to $t$.
\end{theorem}
As for Theorem \ref{th:convUnifOOL}, this results directly leads to Theorem \ref{th:convUnifGUCB}.
To prove Theorem \ref{th:optGUCB}, we proceed as for Theorem~\ref{th:convTau}: we first consider for every $\epsilon>0$ and for every $l\in\{ \epsilon N,\dots,Q^N_1\}$ the number $\Tt^N(l)$ of steps until, using the Good-UCB algorithm, at most $l$ interesting items remain undiscovered on all experts.
Obviously, on $C^N$, the number $\Tt^N(l)$ is upper-bounded by the first time
\[U^N(l) = \inf\left\{t>0 : \forall i\in\{1,\dots,K\}, u^N_{i,n_{i,t}} \leq \frac{l}{N}\right\}\]
when all the upper-bounds of the missing masses are below $l/N$.
Write 
\[U^N(l) = U^N_1(l) + \dots+U^N_K(l)\;,\]
where $U^N_i(l) = n_{i, U^N(l)}$ denotes the number of requests to expert $i\in\{1\dots,K\}$ up to time $U^N(l)$. Four cases can be distinguished:
\begin{itemize}
 \item either the draw does not belong to $C^N$, which has probability $O(N^{-2})$;
 \item or $U^N(l)\geq N^2$, which is easily shown to have probability $O(N^{-2})$; in fact, after $6N\log(N)$ requests to an expert $i$ the probability that there exists any item $(i,x)\in\Omega^N$ which has not been drawn at least twice is $O(N^{-2})$, and otherwise $u^N_{i, N^{3/2}} = O\left( \sqrt{\log(N)/N} \right)<\epsilon$ for $N$ large enough; 
 \item or $U^N_i(l)\leq N^{2/3}$;
 \item or, as $N^{2/3}\leq U^N_i(l)\leq U^N(l)\leq N^2$, and according to Equation~\ref{eq:defC},
 \begin{align*}
  \underline{R}^N_{i,U^N_i(l)} &\geq u^N_{i,U^N_i(l)} - 2\sqrt{\frac{(2/U^N_i(l)+1/N)^2 U^N_i(l) \log(2KN^4)}{2}} - \frac{1}{U^N_i(l)} \\
  &\geq u^N_{i,U^N_i(l)} - C N^{-2/3}\log(KN) \;,
 \end{align*}
for some positive absolute constant $C$. 
Now, remark that
\begin{itemize}
\item $u^N_{i,U^N_i(l)-1} > l$, because the $(U^N_i(l)-1)$-th request to expert $i$ took place before time $U^N(l)$, and at that time expert $i$ had the highest index;
\item $u^N_{i,U^N_i(l)} \geq u^N_{i,U^N_i(l)-1} - C' N^{-2/3}$ for some absolute constant $C'$, as  $n\geq N^{2/3}$, $u^N_{i,n}-u^N_{i,n-1} \geq -1/n$ 
 and \newline $\sqrt{(2/n+1/N)^2 n \log(2KN^4)/2} - \sqrt{(2/(n-1)+1/N)^2 (n-1) \log(2KN^4) /2} \geq -1/n$ because
 \[\frac{\sqrt{(2/n+1/N)^2 n \log(2KN^4)/2}}{\sqrt{(2/(n-1)+1/N)^2 (n-1) \log(2KN^4) /2}} \geq \sqrt{\frac{n}{n+1}}\geq 1-\frac{1}{n}\;,\]
 and $\sqrt{(2/(n-1)+1/N)^2 (n-1) \log(2KN^4) /2}\leq 1$ for $n\geq KN^{3/5}$.
\end{itemize}
% But in that case $u^N_i(t) \geq l - cN\log(N)/\sqrt{B^N}$ because $u^N_i(t)$ can vary at most of $cN\log(N)/\sqrt{B^N}$ in just one step, and if arm $i$ was chosen it means that $u^N_i(t)$ was higher than $l$ the step before (if it had been chosen while $u^N_i(t)<l$, then all upper bounds would have been smaller than $l$ and $U^N(l)$ would be smaller).
Hence $\underline{R}^N_{i,U^N_i(l)}\geq  l - \Delta^N/N$, with $\Delta^N \leq C'' N^{1/3}\log(KN)$ and $C'' = C+C'$. As $\underline{R}^N_k(i) = (Q^N_i - F^N_{i,k})/N$ for every positive integer $k$ , this implies that 
$U^N_i(l) \leq D^N_{i, Q^N_i-l-\Delta^N}$.
\end{itemize}
% Hence, $R^N_k(i) \geq  l - \Delta^N$, with 
% \[\Delta^N = cN\log(N)/\sqrt{B^N} - 2 c\sqrt{\frac{\log(KN^3\log(N))}{N^{2/3}}} = O\left(N^{2/3}\right)\;,\] which implies that 
% $U^N_i(l)\leq T^N_{i, Q^N_i-l-\Delta^N}$ and thus

Altogether, we obtain that on $C^N$:
\[U^N(l) \leq K N^{2/3} + \sum_{i:Q^N_i>l-\Delta^N} D^N_{i,Q^N_i-l-\Delta^N} = K N^{2/3} + T^N(l-\Delta^N)\;.\]
Thus, if $l^N$ is a sequence of integers of $[\epsilon N, Q^N_1]$ such that $l^N/N\to\lambda\in[\epsilon, q_1]$ as $N$ goes to infinity, then
\[\limsup_{N\to\infty} \frac{U^N(l^N)}{N} \leq \lim_{N\to\infty} \frac{T^N(l^N-\Delta^N)}{N} = T(\lambda)\]
according to Lemma~\ref{lem:convT} almost surely, using the Borel-Cantelli lemma and the fact that $\sum_N P(C^N \cup \{U^N(l)\geq N^2\})<\infty$.
This is sufficient to show that $\Tt(l^N)/N$ converges a.s. to $T(\lambda)$ as $N$ goes to infinity.
As for Corollary~\ref{cor:convF}, this implies that $\Ft^N(t^N)/N$ converges a.s. to  $F(t)$ when  $t^N/N$ converges to~$t$.

\section{The open-loop oracle policy}\label{sec:OOL}
An open-loop policy must choose, for each horizon $t$, the respective numbers of requests $(n^N_1,\dots,n^N_K)$ for each distribution (so that $n^N_1+\dots+n^N_K = t^N$) in advance.
It appears here that, in the limit, the \emph{oracle open-loop} (OOL) policy, which makes use of the parameters $(Q^N_1,\dots,Q^N_K)$, is as good as the OCL policy.

Recall that $\underline{R}^N_{i,n^N_i}= (Q^N_i - F^N_{i}(n_i^N))/N$ is the proportion of interesting items not yet found with expert $i$ after $n^N_i$ requests.
Suppose that $t^N/N\to t$, and that $n^N_i/N\to \nu_i$ as $N$ goes to infinity; it is easily shown that almost surely
\[ \lim_{N\to\infty} \underline{R}^N_{i,n_i^N} =\lim_{N\to\infty} \E\left[\underline{R}^N_{i,n_i^N} \right] = \lim_{N\to\infty} \frac{ Q^N_i\left(1-\frac{1}{N}\right)^{n^N_i}}{N} = q_i\exp(-\nu_i)\;.\]
Hence, the proportion of interesting items found with the allocation $(n^N_1,\dots,n^N_K)$ almost surely converges to $\sum_{i=1}^K q_i\left( 1-\exp(-\nu_i) \right)$.
Defining 
\[r(\nu)= \sum_{i=1}^K q_i\exp(-\nu_i)\;,\] it follows that finding the best macroscopic allocation reduces to the following constrained convex minimization problem:
\[\min_{\nu\in\R^K} r(\nu) \hbox{\quad such that } \nu_1+\dots+\nu_K = t \hbox{ and } \forall i,\,\nu_i\geq 0 \;.\]
The solution $r^*(t)$, reached at $\nu=\nu^*(t)$, is easily derived by classical optimization techniques:
\begin{proposition}\label{prop:roo}
For every $i\in\{1,\dots,K\}$, let 
$\uq_i = \exp\left(1/i\times\sum_{k=1}^i \log q_k \right)$ denotes the geometric mean of $q_1,\dots,q_{i}$.
\begin{enumerate}
\item There exists $I(t)\in\{1,\dots, K\}$ such that 
\[\begin{cases}
\forall i\leq I(t),& \nu_i^*(t) =   \frac{t}{I(t)} + \log\frac{q_i}{\uq_{I(t)}} \\
\forall i> I(t),& \nu_i^*(t) = 0 \;.
\end{cases}\]
Hence,
\[r^*(t) = I(t) \uq_{I(t)} \exp\left( -\frac{t}{I(t)} \right) + \sum_{i>I(t)} q_i\;.\]
\item There exists $1=t_1\leq \dots \leq t_K<+\infty$ such that \[\forall t\in[ t_{i},  t_{i+1}[, \;I(t)=i\;.\] 
The $(t_k)_k$ are such that
\[q_{i} + (i-1)\uq_{i-1}\exp\left(-\frac{t_{i}}{i-1}\right) = i\uq_{i} \exp\left(-\frac{t_i}{i}\right)\;.\]
For instance, $t_1 = \log(q_1/q_2)$.
\end{enumerate}
\end{proposition}
\textbf{Proof:}
Introduce the Lagrangian:
\[L(\nu_1,\dots,\nu_K,\lambda, \mu_1,\dots, \mu_K) = \sum_{i=1}^K q_i\exp\left( -\frac{\nu_i}{N} \right) + \lambda \left( \sum_{i=1}^K \nu_i \right) -\sum_{i=1}^K \mu_i \nu_i\]
We need to find the solution of:
\begin{align*}
\forall i\in\{1,\dots, M\},&\;\; -q_i\exp\left( -\nu_i \right) +\lambda - \mu_i = 0\\
&\;\;\sum_{i=1}^K \nu_i = t\\
\forall i\in\{1,\dots,M\},& \;\;\mu_i \nu_i = 0 \hbox{ and } \mu_i \geq 0
\end{align*}
We first obtain that 
\[\nu_i =\log q_i -\log(\lambda - \mu_i)\]
Denoting $A=\{i : \nu_i>0\}$, and using that $i\in A\implies \mu_i=0$ we get:
\[t = \sum_{i\in A} \log(q_i) - |A| \log(\lambda)\]
from which we get
\[-\log(\lambda) = \frac{t}{|A|} - \frac{1}{|A|}\sum_{i\in A}\log q_i, \]
and then for all $i\in A$:
\[\nu_i = \log q_i +\frac{t}{|A|} - \frac{1}{|A|}\sum_{i\in A}\log q_i\;.\]
Next, observe that $\nu_i=0 \iff q_i > \lambda$: in fact, if $\nu_i=0$ then the first equation gives $-q_i +\lambda -\mu_i = 0$, and $0\leq \mu_i = \lambda-q_i$.
Conversely, if $\nu_i>0$ then $\mu_i=0$ and $\nu_i = \log(q_i/\lambda)>0$ implies $q_i> \lambda$.
Thus, there exists $I(t)$ such that $A=\{1,\dots, I(t)\}$, and for all $i\leq I(t)$, 
\[\nu_i = \log\frac{q_i}{\uq_{I(t)}} + \frac{t}{I(t)}\;.\]
Moreover, 
\begin{align*}
r^*(t)  &= r\left(\nu_1,\dots, \nu_{I(t)}, 0,\dots, 0\right) \\
 &= \sum_{i\leq I(t)} q_i \exp\left[-\left(\log\frac{q_i}{\uq_{I(t)}} + \frac{t}{I(t)}\right)\right] +\sum_{i>I(t)} q_i\\
 &= I(t) \uq_{I(t)} \exp\left( -\frac{t}{I(t)} \right)+\sum_{i>I(t)} q_i\;.
\end{align*}
%
%The fact that $I(t)$ is non-decreasing can be seen as follows: if $t<t'$ and if $I(t')<i \leq I(t)$ then $\uq_{I(t')} \geq \uq_{I(t)}$ from which we get $\lambda'>\lambda$, but $q_i<\lambda'$ as $i>I(t')$, and $q_i\leq \lambda$ as $i\leq I(t)$, impossible. WARNING THIS IS FALSE !! important !
%
The instants $(t_i)_{1\leq i\leq K}$ are such that 
\[ (i-1) \uq_{i-1} \exp\left( -\frac{t_i}{i-1} \right)+\sum_{k>i-1} q_k =  i \uq_{i} \exp\left( -\frac{t_i}{i} \right)+\sum_{k>i} q_k\;, \]
which is equivalent to
\[q_i + (i-1)\uq_{i-1}\exp\left(-\frac{t_i}{i-1}\right) = i\uq_i\exp\left(-\frac{t_i}{i}\right)\;.\]
For $i=2$, this gives
\[0 = q_2 + q_1\exp(-\nu_2) - 2\sqrt{q_1q_2} \exp\left( -\frac{\nu_2}{2} \right) = \left( \sqrt{q_2} -\sqrt{q_1\exp\left( -\nu_2 \right)} \right)^2\;,\]
which leads to  $t_1 = \log(q_1/q_2)$.

\begin{theorem}\label{th:OOLoptimal}
The proportion of items found by the open-loop oracle policy uniformly converges to $F$ in the macroscopic limit. 
\end{theorem}
The proportion of interesting items found by the OOL policy is 
\[q - r^*(t) = \sum_{i\leq I(t)} \left[ q_i - \uq_{I(t)}\exp\left( -\frac{t}{I(t)} \right)  \right]
=\sum_{i=1}^K \left( q_i - \Lambda(t) \right)_+
\;,\]
where $\Lambda(t) = \uq_{I(t)}\exp\left( -\frac{t}{I(t)} \right)\in [0, q_{I(t)}]$.
To conclude, it remains only to remark that $\Lambda = T^{-1}$~: in fact, if $\lambda$ is such that $q_{i_0+1}< \lambda\leq q_{i_0}$, then
$I(T(\lambda)) = i_0$ and 
\begin{equation*}
\Lambda\left( T(\lambda) \right) = \uq_{i_0}\exp\left( -\frac{T(\lambda)}{i_0}\right) \\
 = \exp\left( \frac{1}{i_0}\sum_{i\leq i_0} \log q_i \right) \exp\left( -\frac{\sum_{i\leq i_0} \log(q_i / \lambda)}{i_0}\right) = \lambda\;.
\end{equation*}
If $\lambda<q_K$, the same holds with $i_0=K$.

\begin{remark} \label{rmk:1}
Proposition~\ref{prop:roo} permits to compare easily the macroscopic performance of the Good-UCB algorithm with balanced sampling: when all distributions are sampled equally often, the proportion of unseen  interesting items at time $t$ is $\sum_{i=1}^K q_i\exp(-t/K) = K\bar{q}_K\exp(-t/K)$, with $\bar{q}_K = (\sum_{i=1}^K q_i)/K$.
On the other side, for the Good-UCB algorithm, this proportion is 
$I(t) \uq_{I(t)} \exp\left( -t/I(t) \right)$. 
The ratio is a decreasing function of time, but even when $t>t_K$, it takes the nice expression $\bar{q}_K / \uq_K\geq 1$, the ratio of an arithmetic mean with a geometric mean, which is (as expected) as high as the $(q_i)_i$ are unbalanced. 
\end{remark}
\section{Simulations}\label{sec:simus}
\begin{figure}

 \caption{Number of items found by Good-UCB (solid), the oracle (dashed), and uniform sampling (dotted) as a function of time for sizes $N=128, N=500, N=1000$ and $N=10000$ in a $7$-experts setting.}
  \label{fig:nbFound}
  \includegraphics[width=8cm]{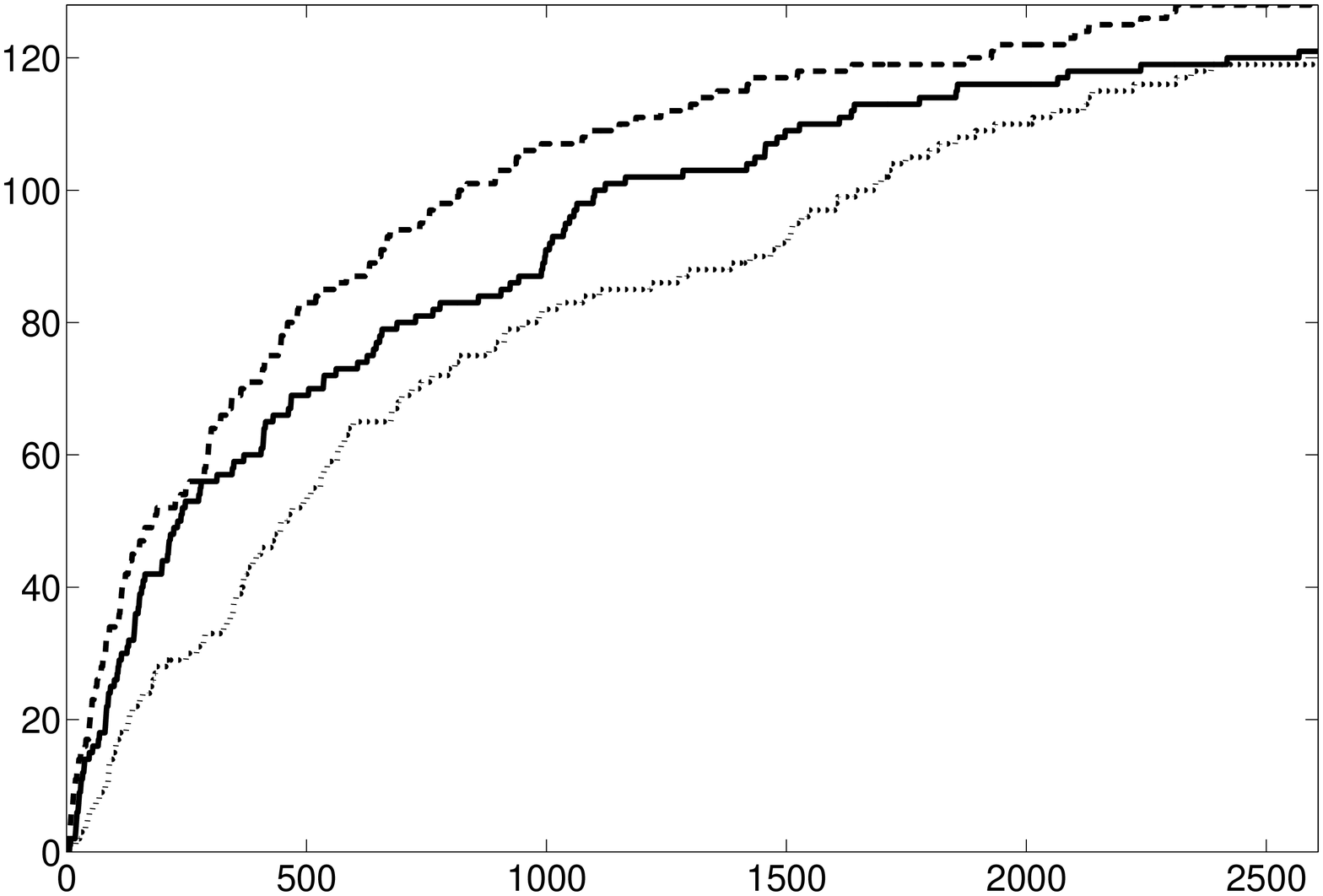}
  \includegraphics[width=8cm]{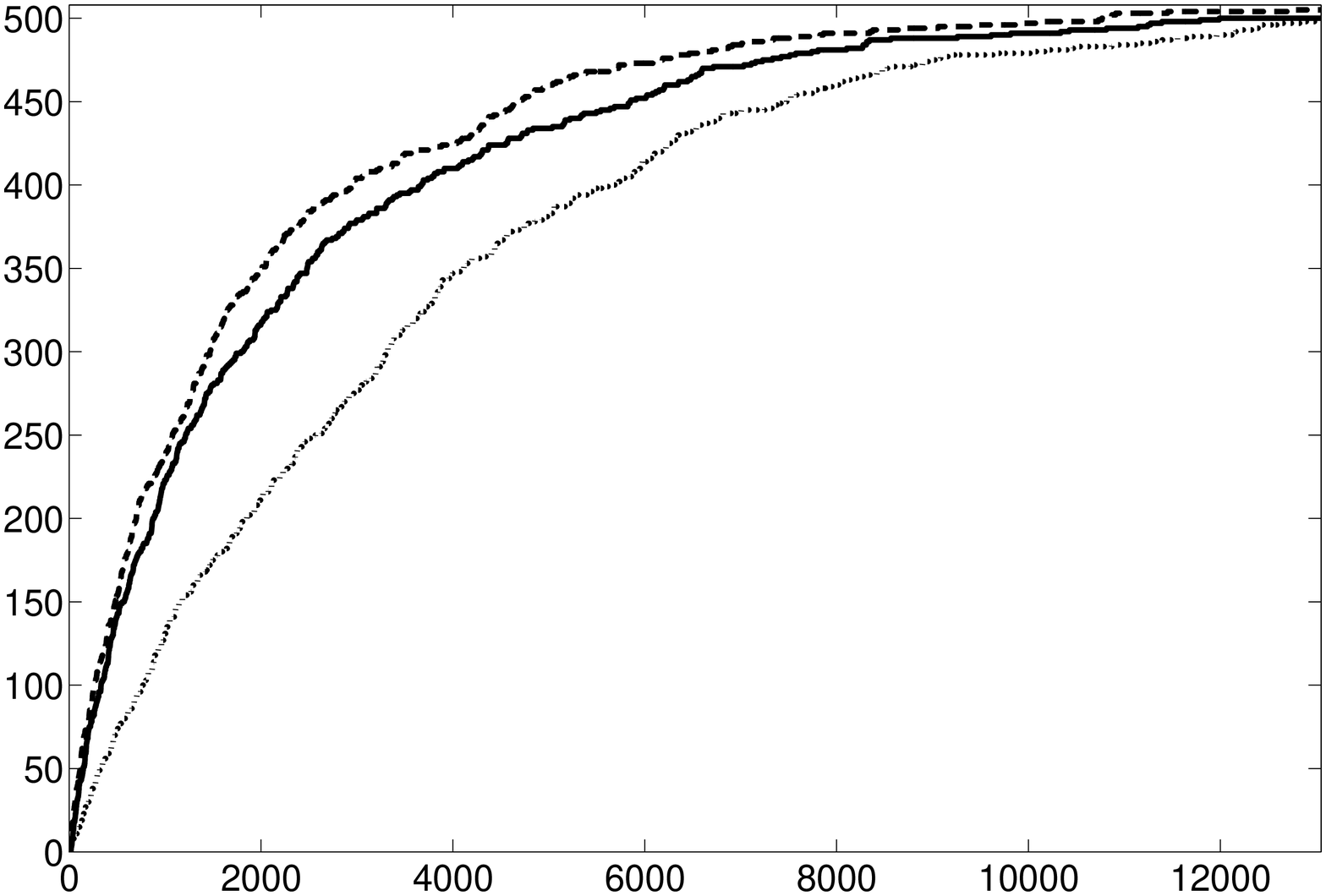} 
 
  \includegraphics[width=8cm]{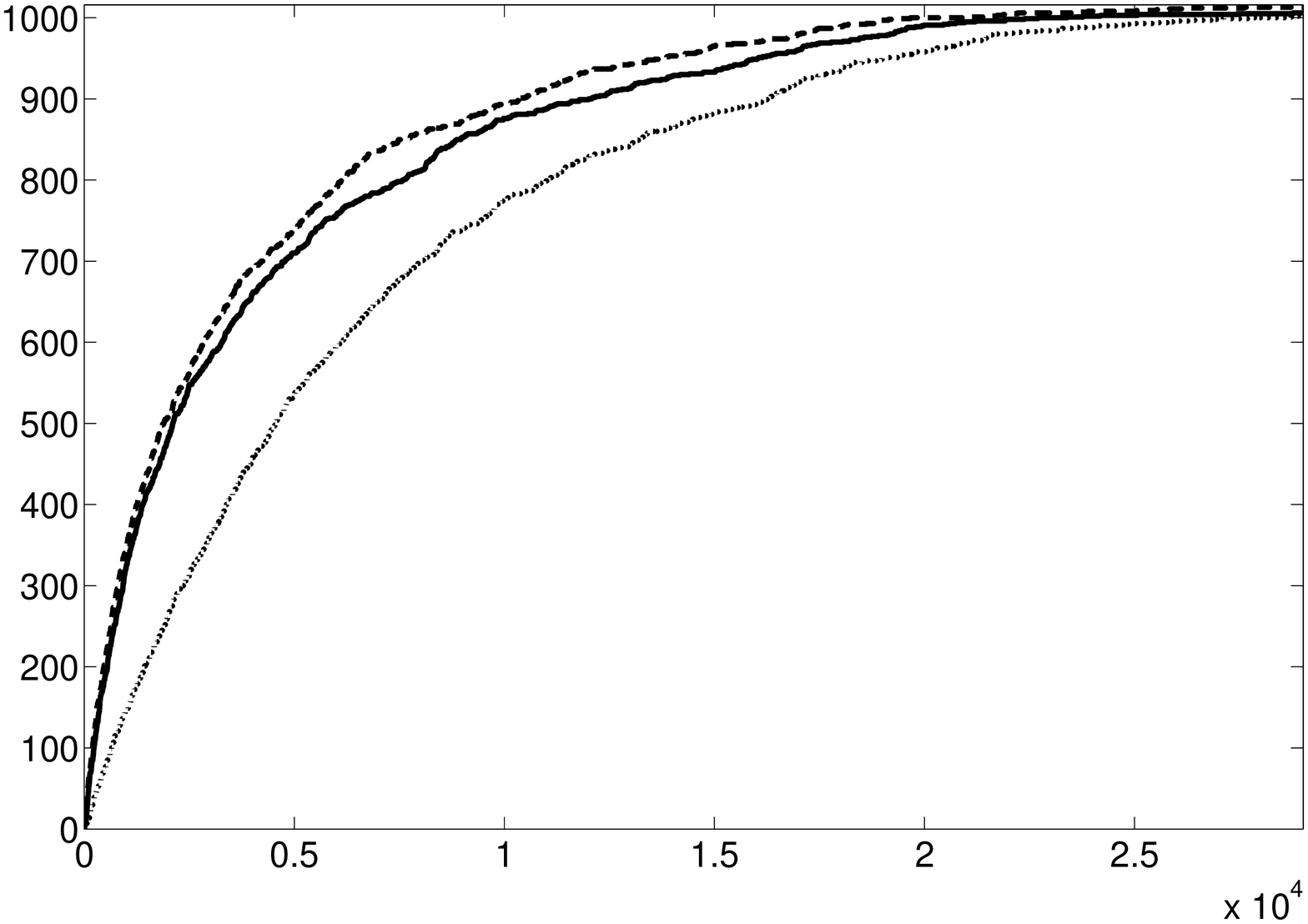}
  \includegraphics[width=8cm]{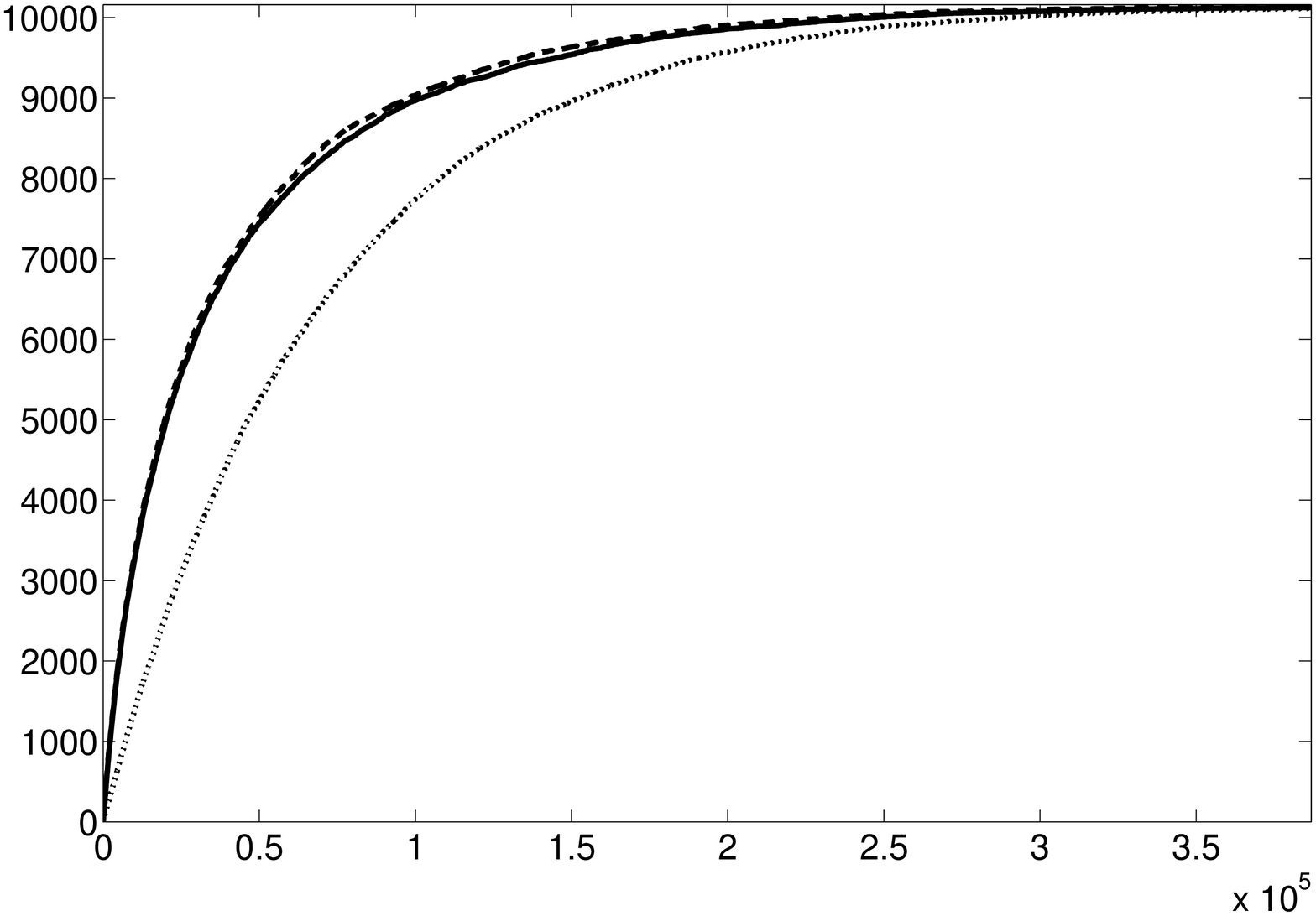}
\end{figure}

We provide a few simulations illustrating the behaviour of the Good-UCB algorithm in practice. 
In order to illustrate the convergence properties shown in Sections~\ref{sec:oracle} and~\ref{sec:goodOptimal}, we first consider an example with $K=7$ different sampling distributions satisfying assumptions [(i),(ii),(iii)], with respective proportions of interesting items $q_1 = 51.2\%, q_2 = 25.6\%, q_3 = 12.8\%, q_4 = 6.4\%, q_5 = 3.2\%, q_6 = 1.6\%$ and $ q_7 = 0.8\%$.
Figure~\ref{fig:nbFound} displays the number of items found as a function of time by the Good-UCB (solid), the oracle (dashed) and a balanced sampling scheme simply alternating between experts (dotted). The results are presented for sizes $N=128, N=500, N=1000$ and $N=10000$, each time for one representative run (averaging over different runs removes the interesting variability of the process).
The convergences proved in Corollary~\ref{cor:convF} and Theorem~\ref{th:optGUCB} are obvious. Moreover, it can be seen that, even for very moderate values of $N$, the Good-UCB significantly outperforms uniform sampling even if it is clearly distanced by the oracle,

For these simulations, the parameter $c$ of Algorithm Good-UCB  has been taken equal to $1/2$, which is a rather conservative choice. In fact, it appears that during all rounds of all runs, all upper-confidence bounds did contain the actual missing mass. Of course, a bolder choice of $c$ can only improve the performance of the algorithm, as long as the confidence level remains sufficient.

%Seb: I removed the following two lines.
%As expected by the completely non-parametric nature of the Good-Turing estimator, experiments show that even for non-uniform %distributions, Good-UCB manages to track the path of the best possible draws.

\begin{figure}
 \caption{Number of prime numbers found by Good-UCB (solid), the oracle (dashed), and uniform sampling (dotted) as a function of time, using geometric experts with means $100, 300, 500, 700$ and $900$, for $c=0.1$ (left) and $c=0.02$ (right).}
  \label{fig:primeNumbers}
  \includegraphics[width=8cm]{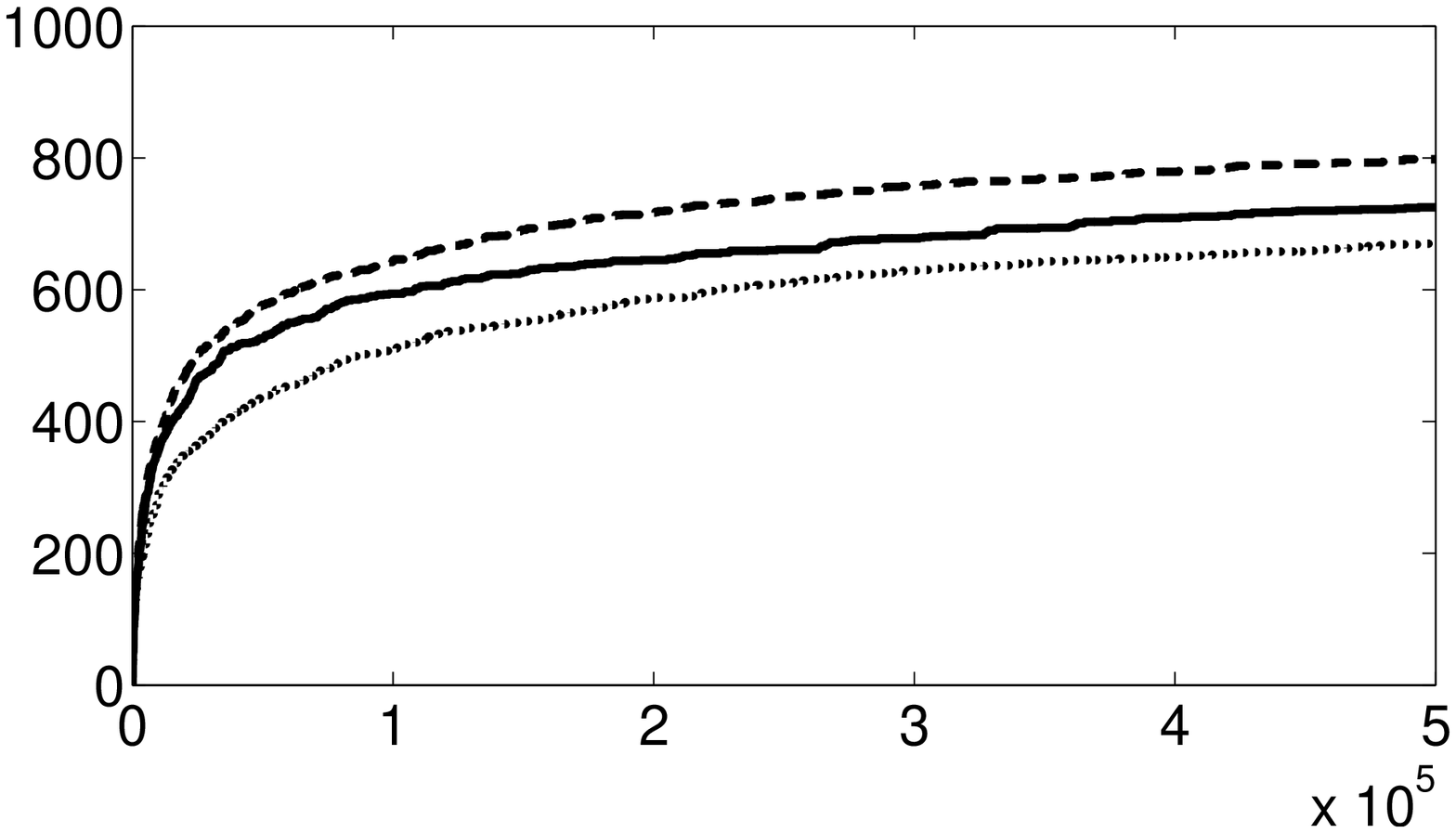}
  \includegraphics[width=8cm]{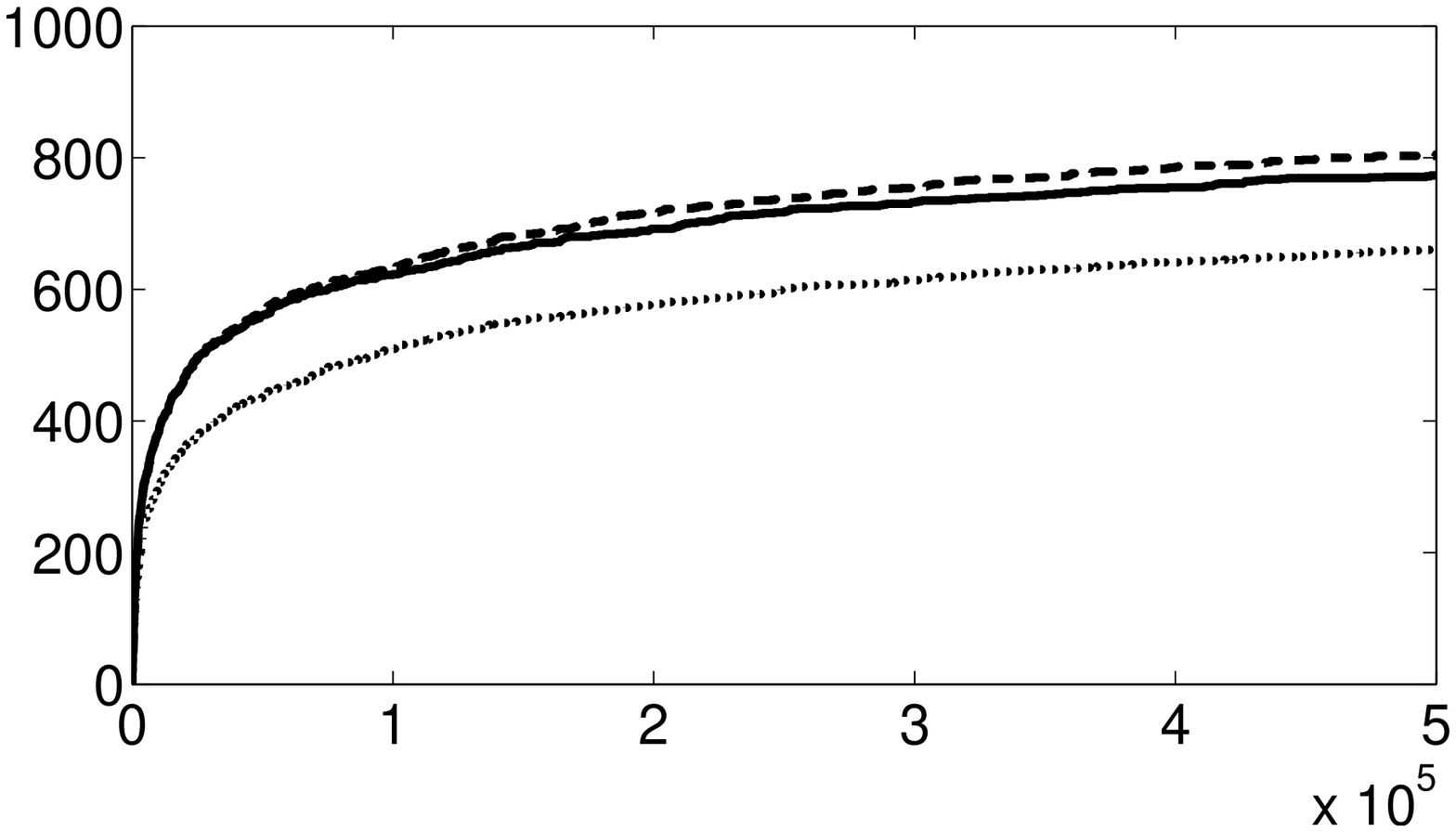}
\end{figure}

In order to illustrate the efficiency of the Good-UCB algorithm in a more difficult setting, which does not satisfy any of the assumptions (i), (ii) and (iii), we also considered the following (artificial) example: $K=5$ probabilistic experts draw independent sequences of geometrically distributed random variables, with expectations $100$, $300$, $500$, $700$ and $900$ respectively. The set of interesting items is the set of prime numbers.
%His goal is to find rapidely as many prime numbers as possible. 
We compare the oracle policy, Good-UCB and uniform  sampling.
The results are displayed in Figure~\ref{fig:primeNumbers}. Even if the difference remains significant between Good-UCB and the oracle, the former still performs significantly better that uniform sampling during the entire discovery process. 
In this example, choosing a smaller parameter $c$ seems to be preferable; this is due to the fact that the proportion of interesting items on each arm is low; in that case, one can show by using tighter concentration inequalities that the concentration of the Good-Turing estimator is actually better than suggested by Proposition~\ref{prop:inegGT}. In fact, this experiment suggests that the value of $c$ should be chosen smaller when the remaining missing mass is small. 

\section{Conclusions} \label{section:conclusions}

This paper introduced an original problem, optimal discovery with probabilistic expert advice. We proposed an algorithm to solve this problem, and showed both analytically and through simulations its efficiency.
%We insist that hypotheses (ii) and (iii) are absolutely not used by the Good-UCB algorithm: they are used only to (considerably)  simplify %the analysis of the oracle policy, and hence to prove the optimality of our algorithm.

This work can be extended along several directions. First, it would be interesting to analyze the behaviour of Good-UCB under less restrictive assumptions on the experts. Note that assumptions (ii) and (iii) are used only to (considerably) simplify the analysis of the oracle policy, and hence to prove the optimality of Good-UCB. Removing assumption (ii) seems fairly straightforward up to the addition of another level of notations. Because of the intricate behavior of oracle policies in that case, it is not clear how the analysis could be carried out if assumption (iii) were to be removed, though it seems reasonable to assume that Good-UCB will still be macroscopically optimal.
% Assumption (i) seems to be the most challenging to remove. 
Good-UCB is designed to work even when assumption (i) is not satisfied, but the analysis is complicated by the explicit dependency between the missing mass of the different experts.

One may also wonder whether it would be possible to obtain optimal {\it rates of convergence} (in the macroscopic limit sense) for this problem, and whether Good-UCB is optimal in that sense too. 
Finally, another macroscopic limit deserves to be investigated, where the number of interesting items for each arm remains constant, while $N$ and $n$ go to infinity. Note that  in such a case,  a Poisson regime appears. The analysis of Good-UCB might be possible by using a better concentration bound for the Good-Turing estimator such as the Boucheron-Massart-Lugosi inequality~\cite{BoucheronLugosiMassart09selfBounding}. This could also contribute to explain why, in the second experiment presented in Section~\ref{sec:simus}, the parameter $c$ should be chosen decreasing with time.

\section{Appendix}
\begin{lemma}\label{lem:harmo}
For all $1\leq k\leq n$, 
\begin{equation}\label{eq:harmonic}
-\frac{1}{k}+\log\frac{n}{k}  \leq \sum_{j=k+1}^n \frac{1}{j} \leq \log\frac{n}{k}
\end{equation}
\end{lemma}
\noindent\textbf{Proof:}
The standard sum/integral comparison yields
\[\log\frac{n+1}{k+1}\leq\sum_{j=k+1}^n \frac{1}{j} \leq \log\frac{n}{k}\, \]
but 
\[\log\frac{n+1}{k+1} = \log\frac{n}{k} + \log\left( 1+\frac{1}{n+1} \right) - \log\left( 1+\frac{1}{k+1} \right)  \geq \log\frac{n}{k} + 0 -\frac{1}{k}\;. \]

\vspace{1cm}\noindent\textbf{Proof of Lemma~\ref{lem:convFL}:}
For $(i)$, it suffices to notice that 
\[\frac{1}{N}G^N(l^N) = \sum_{i=1}^K \left( \frac{Q^N_i}{N} - \frac{l^N}{N}  \right)_+ \to \sum_{i=1}^K \left( q_i - \lambda  \right)_+ \]
if $l^N/N\to\lambda$.
Moreover, the same argument shows that 
\[\limsup \frac{1}{N}G^N(l^N)  = G\left(\limsup  \frac{l_N}{N} \right) \hbox{ and } \liminf \frac{1}{N}G^N(l^N)  = G\left(\liminf  \frac{l_N}{N} \right)\;,\]
Hence, if $G^N(l^N)/N$ converges, the limit belongs to $[0,q]$ and thus can be written $G(\lambda)$ for some $\lambda\in[0,q_1]$.
Thus, $G\left(\liminf  \frac{l_N}{N} \right)  = G \left(\limsup  \frac{l_N}{N} \right) = G(\lambda)$, which implies that $l^N/N\to \lambda$ as $G$ a continuous bijection.

Concerning $(ii)$: if $f^N/N\to \phi$, then the fact that $|G^N(L^N(f^N)) - f^N|<K$ implies that that $G^N(L^N(f^N)) /N \to \phi$ and thus, by $(i)$, that $L^N(f^N)/N$ converges to a value $\lambda$ such that $G(\lambda)=\phi$, i.e. $\lambda = L(\phi)$. The reciprocal (which is not used in the sequel) is left to the reader.

\vspace{1cm}\noindent\textbf{Proof of Lemma~\ref{lem:convT}:} Since 
\[T^N(l^N) = \sum_{i: Q^N_i>l^N} D^N_{i, Q^N_i-l^N}\;,\]
it suffices to show that for every expert $i\in\{1,\dots,K\}$,  $D^N_{i, Q^N_i-l^N}/N$ converges almost surely to $\log(q_i/\lambda)$ as $N$ goes to infinity.
Write
\begin{equation}\label{eq:WsumS}
W^N_{i,l^N} = \frac{1}{N}\left(D^N_{i, Q^N_i-l^N} - \E\left[D^N_{i, Q^N_i-l^N}\right]\right) = 
\frac{1}{N}\sum_{k=1}^{Q^N_i-l^N-1} \left(S^N_{i, k} - \E\left[S^N_{i, k}\right]\right) \;.\end{equation}
For every positive integer $d$ and for $k\in\{1,\dots, l^N-1\}$, elementary manipulations of the geometric distribution yield that \[\E\left[\left(S^N_{i, k} - \E\left[S^N_{i, k} \right]\right)^d\right] \leq \E\left[\left(S^N_{i, l^N} - \E\left[S^N_{i, l^N} \right]\right)^d\right] \leq \frac{c(d)}{(l^N/N)^d} \leq \frac{2c(d)}{\lambda^4}\]
for some positive constant $c(d)$ depending only on $d$, and for $N$ large enough. 
Hence, taking~\eqref{eq:WsumS} to the fourth power and developing yields
\[\E\left[\left(W^N_{i,l^N}\right)^4\right] \leq \frac{c'}{N^2\lambda^4}\] for some positive constant $c'$. Using Markov's inequality together with the Borel-Cantelli lemma, this permits to show that $W^N_{i,l^N}$ converges almost surely to $0$ as $N$ goes to infinity.
But
\[ \frac{1}{N}\E\left[D^N_{i, Q^N_i-l^N}\right] = \frac{1}{Q^N_1} + \dots +\frac{1}{l^N+1} = \log\frac{Q^N_i}{l^N} - \epsilon^N\;, \]
with $0\leq \epsilon^N \leq 1/l^N$ according to Lemma~\ref{lem:harmo}, and thus 
\[\frac{1}{N}\E\left[D^N_{i, Q^N_i-l^N}\right]\to \lim_{N\to\infty} \log\left( \frac{Q^N_i/N}{l^N/N} \right) = \log(q_i/\lambda)\;,\]
which concludes the proof.

\bibliographystyle{alpha}
\bibliography{biblio}

\newcommand{\etalchar}[1]{$^{#1}$}
\begin{thebibliography}{FBED{\etalchar{+}}10}

\bibitem[ACBF02]{AuerEtAl02FiniteTime}
P.~Auer, N.~Cesa-Bianchi, and P.~Fischer.
\newblock {Finite-time analysis of the multiarmed bandit problem}.
\newblock {\em Machine Learning}, 47(2):235--256, 2002.

\bibitem[Agr95]{Agrawal:95}
R.~Agrawal.
\newblock {Sample mean based index policies with O(log n) regret for the
  multi-armed bandit problem}.
\newblock {\em Advances in Applied Probability}, 27(4):1054--1078, 1995.

\bibitem[BLM09]{BoucheronLugosiMassart09selfBounding}
S.~Boucheron, G.~Lugosi, and P.~Massart.
\newblock On concentration of self-bounding functions.
\newblock {\em Electron. J. Probab.}, 14:no. 64, 1884--1899, 2009.

\bibitem[FB11]{Fonteneau11phd}
F.~Fonteneau-Belmudes.
\newblock {\em Identification of dangerous contingencies for large scale power
  system security assessment}.
\newblock PhD thesis, University of Li\`ege, 2011.

\bibitem[FBED{\etalchar{+}}10]{FonteneauErnstDruetPanciaticiWehenkel10power}
F.~Fonteneau-Belmudes, D.~Ernst, C.~Druet, P.~Panciatici, and L.~Wehenkel.
\newblock Consequence driven decomposition of large-scale power system security
  analysis.
\newblock In {\em Proceedings of the 2010 IREP Symposium - Bulk Power Systems
  Dynamics and Control - VIII}, Buzios, Rio de Janeiro, Brazil, August 2010.

\bibitem[Goo53]{Good53GT}
I.J. Good.
\newblock The population frequencies of species and the estimation of
  population parameters.
\newblock {\em Biometrika}, 40:237--264, 1953.

\bibitem[GS95]{Gale95GT}
W.A. Gale and G.~Sampson.
\newblock Good-turing frequency estimation without tears.
\newblock {\em Journal of Quantitative Linguistics}, 2(3):217--237, 1995.

\bibitem[McD89]{Mcdiarmid89boundedDiff}
C.~McDiarmid.
\newblock On the method of bounded differences.
\newblock In {\em Surveys in combinatorics, 1989 ({N}orwich, 1989)}, volume 141
  of {\em London Math. Soc. Lecture Note Ser.}, pages 148--188. Cambridge Univ.
  Press, Cambridge, 1989.

\bibitem[MS00]{McallesterSchapire00GoodTuring}
D.A. McAllester and R.E. Schapire.
\newblock On the convergence rate of good-turing estimators.
\newblock In {\em COLT}, pages 1--6, 2000.

\bibitem[OSZ03]{orlitsky2003always}
A.~Orlitsky, N.P. Santhanam, and J.~Zhang.
\newblock {Always Good Turing: asymptotically optimal probability estimation}.
\newblock In {\em FOCS '03: Proceedings of the 44th Annual IEEE Symposium on
  Foundations of Computer Science}, pages 179+, Washington, DC, USA, 2003. IEEE
  Computer Society.

\end{thebibliography}
\end{document}